\documentclass[reqno]{amsart}
\usepackage{amsfonts}
\usepackage{amssymb}
\usepackage{amsthm}
\usepackage{upref}
\makeatletter
\def\LaTeX{\leavevmode L\raise.42ex
    \hbox{\kern-.3em\size{\sf@size}{0pt}\selectfont A}\kern-.15em\TeX}
 \makeatother
 
\numberwithin{equation}{section}

\newtheorem{lemma}{Lemma}[section]
\newtheorem{theorem}[lemma]{Theorem} 

\newtheorem{proposition}[lemma]{Proposition}

\theoremstyle{definition}

\newcommand{\HH}{\mathsf{H}}
 
  \newcommand{\e}{\eqref}

\newcommand{\q}{\quad}

\newcommand{\wt}{\widetilde}

\renewcommand\Re{\operatorname{Re}}

\def\qqq{\mathrel{\subset\mkern-15mu\lower.38ex\hbox{${\scriptscriptstyle\rightarrow}$}}}
\let\cal\mathcal

\let\Bbb\mathbb

\begin{document}

\title[Diagonalization   of Hankel operators]
{A commutator method for the diagonalization \\ of Hankel operators}
\author{ D. R. Yafaev}
\address{ IRMAR, Universit\'{e} de Rennes I\\ Campus de
  Beaulieu, 35042 Rennes Cedex, FRANCE}
\email{yafaev@univ-rennes1.fr}

\dedicatory{To the memory of Mikhail Shl\"emovich Birman}

\begin{abstract}
We present a method for the explicit diagonalization of some Hankel operators. This method allows us to recover classical results on the diagonalization of Hankel operators with the absolutely continuous spectrum. It leads also to new results. Our approach relies on the commutation of a Hankel operator  with some differential operator   of   second order.
  \end{abstract}

\maketitle

% \thispagestyle{empty}

%************************************************************
\section{Introduction}  
%***********************************************************

{\bf 1.1.}
Hankel operators can be defined (see, e.g., book \cite{Pe}) as integral operators in the space $L^2({\Bbb R}_{+})$ whose kernels depend on the sum of variables only. Thus, a Hankel operator $A$  is defined by the formula
\begin{equation}
(A f)(x)= \int_{0}^\infty a(x+y) f(y) dy.
 \label{eq:H1}\end{equation}
Of course, $A$  is self-adjoint if $a=\bar{a}$.  If 
 \[
 \int _0^\infty |a(x)|^2 x dx< \infty,
 \]
 then $A$ belongs to the Hilbert-Schmidt class. This condition is   satisfied if, for example, the function $a$ is continuous, it is not too singular at $x=0$ and decays sufficiently rapidly as $x\to\infty$.
  On the contrary, if $a(x)\sim a_{0}x^{-1}$ as $x\to 0$ or (and) $a(x)\sim a_{\infty}x^{-1}$ as $x\to\infty$, then the operator $A$ is no longer compact although it remains bounded. A general philosophy (see paper \cite{Ho} by
J.~S.~Howland) is that each of these singularities gives rise to the branch $[0,a_{0} \pi]$ or (and) $[0,a_\infty \pi]$ of the simple absolutely continuous spectrum.
 
 There are very few examples where the operator $A$ can be explicitly diagonalized, that  is its exact eigenfunctions can be found. The first result is due to F. Mehler \cite{Me} who considered the case $a(x)=(x+2)^{-1}$. He has shown that   functions
  \begin{equation}
  \psi _k(x)= \big(  k \tanh   \pi  k \big)^{1/2} P_{-1/2+i k}(x+1), \q \lambda=  \pi/ \cosh \pi k, \q k \ > 0,
 \label{eq:M1}\end{equation}
where $P_{-1/2+i k}$ is the Legendre function (see \cite{BE}, Chapter~3), satisfy  equations $A\psi_{k}=\lambda \psi_k$. The functions $\psi_{k}$ are usually parametrized by the quasimomentum $k$ related to $\lambda=\lambda(k)$ by formula \e{eq:M1}.
The      operator  $U: L^2 ({\Bbb R}_{+}) \to L^2 ({\Bbb R}_{+})$ defined\footnote{A precise definition of the operator $U$ can be given in terms of the corresponding sesquilinear form. }  by the equality
\begin{equation}
(U f) (k)= \int_{0}^\infty  {  \psi}_k (x)  f(x) dx
 \label{eq:U}\end{equation}
 is unitary. Observe that $\lambda(k)$ is a one-to-one mapping of ${\Bbb R}_{+}$ on $(0,\pi)$ and  that $(ULf )(k)=\lambda (k) f(k)$ which implies
   that the spectrum of the operator $A$ is simple, absolutely continuous and coincides with the interval $[0, \pi]$.  
   
   Below we use the term  ``eigenfunction" for $  \psi _k$ (although $  \psi _k\not\in L^2 ({\Bbb R}_{+})$) such that $A\psi_{k}=\lambda \psi_k$    for the   spectral parameter $\lambda$ from the continuous spectrum of the operator $A$. 
   {\it By definition}, we   also say that   
eigenfunctions   
$\psi_k$ of the continuous spectrum   are orthogonal, normalized and the  set of all $\psi_{k}$ is complete if the corresponding operator \e{eq:U} is unitary
(if $A$    has no point spectrum).

The next result is due to  W. Magnus \cite{Ma} who considered the case $a(x)=x^{-1} e^{-x/2}$. A more general result of the same type was obtained by M. Rosenblum \cite{Ro} who has diagonalized the operator $A$ with kernel
\begin{equation}
a(x)=\Gamma(1+ \beta) x^{-1}W_{-\beta, 1/2}(x), \q \beta\in{\Bbb R}, \q \beta\neq -1,-2,\ldots,
 \label{eq:W1}\end{equation}
 where $W_{-\beta, 1/2}$ is the Whittaker function (see \cite{BE}, Chapter~6) and $\Gamma$ is the gamma function. Note that $W_{0, 1/2}(x)=e^{-x/2}$. The spectrum of the operator $A$ with such kernel is again simple and, up to a finite number of eigenvalues, it is absolutely continuous and coincides with the interval $[0, \pi]$. Its ``normalized eigenfunctions" are expressed in terms of the Whittaker functions
 \begin{equation}
\psi_{k} (x)= (2\pi)^{-1} \sqrt{k   | \Gamma (1/2 - ik +\beta)| \sinh 2\pi k}    x^{-1}W_{-\beta, i k}( x)  ,\q k>0.
 \label{eq:Whi}\end{equation}

 Observe that the function $a(x)=(x+2)^{-1}$ is singular at $x=\infty$ and   eigenfunctions  \e{eq:M1} decay as linear combinations of $x^{-1/2\pm ik}$ as $x\to \infty$ while function \e{eq:W1} is singular at $x=0$ and   eigenfunctions \e{eq:Whi}  behave as linear combinations of the same functions $x^{-1/2\pm ik}$ as $x\to 0$.
 
  We note also a simple case $a(x)=x^{-1}$ where the operator $A$ is directly diagonalized (see paper \cite{Ca} by T.~Carleman) by the Mellin transform.  In this case the spectrum of $A$   has multiplicity $2$ (because of the singularities of $a(x)$ both at $x=0$ and at $x=\infty$), it is absolutely continuous and coincides with the interval $[0, \pi]$. The eigenfunctions of  the Carleman operator  equal $x^{-1/2\pm ik}$ (up to a normalization).
  
  We emphasize a parallelism of theories of singular differential  operators and Hankel operators with singular kernels. Thus, the functions $x^{-1/2\pm  ik}$ play (both for $x\to \infty$ and $x\to 0$) for Hankel operators the   role of exponential functions $e^{\pm ikx}$ for differential   operators of second order. From this point of view, the Carleman operator plays the role of the operator $-d^2/ dx^2$ in the space $L^2({\Bbb R})$.
 
 \medskip
 
 {\bf 1.2.}
 In the author's opinion,
the reason why in  the cases described above eigenfunctions of a Hankel operator can be found explicitly remained unclarified.  Our approach shows   that    all diagonalizable Hankel operators $A$ commute with  differential operators 
\begin{equation}
L = - \frac{d}{dx}(x^2+\gamma x) \frac{d}{dx} + \alpha x ^2 + \beta  x
 \label{eq:L1}\end{equation}
 for suitably chosen parameters $\alpha\geq 0,\beta\in{\Bbb R}$ and $\gamma\geq 0$. 
 Thus, eigenfunctions of the operators $A$ and $L$ are the same which allows us to diagonalize the operator $A$. 
 
   Hopefully the commutator method will be applied to other kernels $a$. In this paper we use  the commutator method to find in subs.~4.4 eigenfunctions of a new Hankel operator with kernel
 \begin{equation}
a(x)=   \sqrt{\frac{8}{x}} K_{1} (\sqrt{8x}),
 \label{eq:H2}\end{equation}
  where $K_{1} $ is the MacDonald function (see \cite{BE}, Chapter~7). Similarly to \e{eq:W1}, this function decays exponentially as $x\to \infty$ and $a(x)\sim x^{-1}$ as $x\to 0$. An example of a different nature are Hankel operators with regular kernels; such operators are compact.
 
  % We recall that $K_{1}(z)=. H_{1}^{(1)}(z)$ where $H_{1}^{(1)}$ is the Hankel function.
 
 Note that operator \e{eq:L1} for $\gamma=0$ and $\alpha>0$ appeared already in \cite{Ro}. Actually, M.~Rosenblum proceeded from the identity 
\begin{equation}
\Gamma(1+ \beta)\int_{0}^\infty (x+y)^{-1}W_{-\beta, 1/2}(x+y) y^{-1}W_{-\beta, i k}( y) dy
=\frac{\pi}{\cosh \pi k}   x^{-1}W_{-\beta, i k}( x)  
 \label{eq:Sh}\end{equation}
 found earlier by H. Shanker in \cite{Sh}. This identity shows that functions  \e{eq:Whi} are eigenfunctions of the Hankel operator with kernel \e{eq:W1}. M.~Rosenblum
  observed that functions  \e{eq:Whi}  are also eigenfunctions of   operator \e{eq:L1}  for $\gamma=0$ and $\alpha=1/4$. Since eigenfunctions of the self-adjoint differential operator $L$ are orthogonal and complete, the same is true for eigenfunctions of the Hankel operator   $A$ with kernel \e{eq:W1}.
   This yields the diagonalization of this operator. 
  
  Our approach is somewhat different. We prove the relation $LA=AL$ which shows that eigenfunctions of the operators $L$ and $A$ are the same. In particular, we obtain identity \e{eq:Sh} without a recourse to  the theory of special functions.
 
It is well-known   that the integrability   of  differential equations  of second order in terms of special functions has a deep group-theoretical interpretation (see, e.g.,  book \cite{Vi} by N.~Ya.~Vilenkin).  As far as Hankel operators are concerned, it is evident that the diagonalization of the Carleman operator can be explained by its invariance with respect to the group of dilations. The relation $LA=AL$   means that the operator $A$ is invariant with respect to the group $\exp(-itL)$. In contrast to the Carleman operator, for other Hankel operators this invariance does not look obvious.

A commutator scheme is presented in Section~2 while specific examples of kernels singular at $x=\infty$ and $x=0$ are discussed in Sections~3 and 4, respectively. Hankel operators with regular kernels are considered in Section~5.

%%%%%%%%%%%%%%%%%%%%%%%%%%%
\section{Commutator method}  
%%%%%%%%%%%%%%%%%%%%%%%%%%%%

 {\bf 2.1.}
  For a moment, we consider the operator  $L$ defined by formula \e{eq:L1}  as a differential operator on the class
       $  C ^2 ({\Bbb R}_{+})$, but later it will be defined as a self-adjoint operator in the space  $  L^2 ({\Bbb R}_{+})$. Let the operator $A$ be given by formula \e{eq:H1} where $a\in C ^2 ({\Bbb R}_{+})$. 
       
       Let us commute the operators $A$ and $L$. Suppose that $f\in C ^2 ({\Bbb R}_{+})$ and that
        \begin{equation}
 \lim_{y\to 0}  (y^2+ \gamma y) f (y) =   \lim_{y\to 0}  (y^2+ \gamma y) f' (y)=  0
\label{eq:L4} \end{equation}
as well as
     \begin{equation}
 \lim_{y\to\infty}    a'(x+y) (y^2+ \gamma y) f (y) =   \lim_{y\to\infty}  a(x+y) (y^2+ \gamma y)  f' (y)=  0
\label{eq:L4a} \end{equation}
 for all $x\geq 0$.    Then integrating by parts, we find that 
 \[
(AL f)  (x)= \int_{0}^\infty \Big(- \frac{\partial}{\partial y}\big( (y^2+\gamma y) a' (x+y) \big)+ a(x+y)(\alpha y^2+ \beta y) \Big) f(y) dy.
\]
It follows that
 \[
(( LA-AL) f)  (x)= \int_{0}^\infty q (x,y) f(y) dy
\]
where
\begin{align*}
  q (x,y)=& - \frac{\partial}{\partial x}\big( (x^2+\gamma x) a' (x+y) \big) +  \frac{\partial}{\partial y}\big( (y^2+\gamma y) a' (x+y) \big) 
  \\
&  +  (\alpha x^2 -\alpha y^2+ \beta x- \beta y) a(x+y) 
\nonumber\\
=& (x-y) \Big(-(z+\gamma) a''(z) -2a'(z) + (\alpha z +\beta) a(z) \Big)  
 \end{align*}
and $z=x+y$. Thus, we arrive at the following general result.

\begin{theorem}\label{pr1}
Suppose that   kernel $a$ of a Hankel operator $A$ satisfies the differential  equation 
\begin{equation}
 -(x+\gamma) a''(x) -2a'(x) + (\alpha x+\beta) a(x)=0.
\label{eq:E1} \end{equation}
Let $f\in C ^2 ({\Bbb R}_{+})$  and let conditions \e{eq:L4} and \e{eq:L4a} hold. Then 
\begin{equation}
( LA-AL) f=0.
\label{eq:comm} \end{equation}
\end{theorem}

Note that
after a change of variables 
\begin{equation}
  a(x) = (x+\gamma)^{-1} b(x+\gamma)
\label{eq:E2} \end{equation}
in \e{eq:E1},  we get the Schr\"odinger equation with the Coulomb potential
\begin{equation}
 -  b''(r)  + (\alpha  +\beta r^{-1}) b(r)=0.
\label{eq:E3} \end{equation}

\medskip
 
 {\bf 2.2.}
 In specific examples below, we are going to use Theorem~\ref{pr1} in the following way. If $L$ is self-adjoint and has a simple  spectrum, then the equality $LA = AL $ shows that $A$ is a function   ${\sf F}$ of $L$, i.e., the operators $A$ and $L$ have common eigenfunctions. For a calculation of the  function   ${\sf F}$, we argue as follows. Suppose that a function $\psi_{\mu}$ satisfies conditions \e{eq:L4}, \e{eq:L4a} and the equation
\begin{equation}
 -\big((x^2+\gamma x)  \psi_{\mu}'(x)\big)'  + (\alpha x^2+\beta x) \psi_{\mu}(x)=\mu \psi_{\mu}(x).
\label{eq:Y1} \end{equation}
Then according to equality \e{eq:comm} the same equation holds for the function $A  \psi_\mu$ and hence, for some numbers $\lambda=\lambda_{\mu} $ and $\check{\lambda} =\check{\lambda}_{\mu}$, 
\begin{equation}
(A \psi_{\mu})(x)= \lambda\psi_{\mu}(x)+ \check{\lambda} \check{\psi}_{\mu}(x)
\label{eq:A1} \end{equation}  
where $\check{\psi}_{\mu}$ is a solution of the equation $L\check{\psi}_{\mu}=\mu\check{\psi}_{\mu}$  linearly independent of $\psi_{\mu}$. Further,  comparing asymptotics of the functions  $  \psi_{\mu}(x)$, $\check{\psi}_{\mu}(x)$ and $ (A \psi_{\mu})(x)$  as $x\to 0$ and $x\to \infty$, we see that $\check{\lambda}=0$ and find $\lambda={\sf F} (\mu) $ as a function of $\mu$.  
 Finally, if $\psi_{\mu}$ belong to the domain of some self-adjoint realization of the differential operator $L$, then, for a proper normalization of functions $ \psi_{\mu}$,
   the system of all $ \psi_{\mu}$ is orthogonal and complete.  In this case $A={\sf F} (L)$.
   Note that this approach allows one to avoid     precise definitions of   commutators and references to the functional analysis.
 
It turns out that in all our applications   ${\sf F} (\mu)=  \pi / \cosh\big(\pi \sqrt{ \mu -1/4 }\big)$, and hence
\[
A= \pi / \cosh\big(\pi \sqrt{L-1/4 }\big).
\]
  Actually, it is somewhat more convenient  to parametrize  eigenfunctions by the quasimomentum $k > 0$ related to $\mu$ and $\lambda$ by the formulas
\begin{equation}
 \mu=k^2+1/4\in (1/4, \infty),\q \lambda= \pi / \cosh \pi k \in (0,\pi).
\label{eq:P1} \end{equation}
Note that $ \psi_k(x)$, $ \psi_{\mu}(x)$ and $ \psi_\lambda(x)$ denote the same function provided the parameters $k$, $\mu$ and $\lambda$ are related by formulas \e{eq:P1}.

The  operator $U$   defined by formula \e{eq:U}  is unitary and  the operator $UA U^* $ acts in $L^2 ({\Bbb R}_{+})$ as   multiplication by the function $\lambda(k)= \pi / \cosh \pi k$. Indeed,  according to the Fubini theorem it follows   from the equation $A\psi_{k }=\lambda(k)\psi_{k}$ that for  $g\in C_{0}^\infty ({\Bbb R}_{+})$ 
\begin{align*}
(A U^* g) (x)&=\int_{0}^\infty d k g(k) \int_{0}^\infty dy  a(x+y) {  \psi}_k (y)
\\
&= \int_{0}^\infty \lambda (k) {  \psi}_k (x) g(k)d k= (U^* (\lambda  g))(x),
\end{align*}
or equivalently
\begin{equation}
 ( U A f)(k)= \lambda(k) ( U f)(k),\q \forall f\in L^2 ({\Bbb R}_{+}).
\label{eq:CC} \end{equation}
Since $\lambda:{\Bbb R}_{+}\to ( 0,\pi )$ is a  smooth one-to-one mapping, the operator $A$ has the simple absolutely continuous spectrum $[0,\pi]$.

To realize this scheme, it is convenient to study the cases of singularities at $x=\infty$ when $\gamma>0$ and at $x=0$ when $\gamma=0$ separately.

% It follows that
%   the operator $U$   defined by formula \e{eq:U} where ${\pmb \psi}_{\lambda}(x)=F'(\mu)^{-1/2} \psi_{\mu}(x)$ is unitary and diagonalizes the operator $A$. Indeed,  according to the Fubini theorem for $g\in L^2 (0, \pi)$, we have
%(A U^* g) (x)=\int_{0}^\pi d\lambda g(\lambda) \int_{0}^\infty dy  a(x+y) {\pmb \psi}_{\lambda} (y)= \int_{0}^\pi  \lambda  {\pmb \psi}_{\lambda} (x) g(\lambda)d\lambda= (U^* (\lambda g))(x).\]
% In this case $A= F (L)$.

%%%%%%%%%%%%%%%%%%%%%%%%%%%
\section{Singularity at infinity}  
%%%%%%%%%%%%%%%%%%%%%%%%%%%%

{\bf 3.1.}
Set $\gamma=2$.
 We first suppose that $\alpha=\beta=0$. Then the function $a(x)= (x+2)^{-1}$ satisfies equation \e{eq:E1}, and the corresponding operator
\[
L = - \frac{d}{dx} p (  x) \frac{d}{dx} \q {\rm where} \q p(x)=x^2+ 2 x. 
\]

Let $P_\nu(z)$ and $Q_{\nu}(z) $ be the Legendre functions (see, e.g., \cite{BE}, Ch.~3) of the first and second kinds, respectively. They are defined as solutions of the   equation
\[
(1-z^2) u'' (z) -2z u'(z)+ \nu(\nu+1) u(z)=0, \q z>1,
\]
satisfying the conditions $P_\nu(1)=1$ and $Q_\nu(z)=-2^{-1}\ln (z-1)+ c_{\nu}$ as $z\to 1+0$ (the value of the number $c_{\nu}$ is inessential).  Then the functions  
$P_{-1/2+ i k}(x+1)$ and $Q_{-1/2+ ik}(x+1)$ satisfy   the equation $Lu = (k^2+1/4) u $.
We  also note  that (see formulas (2.10.2) and (2.10.5) of \cite{BE})
\begin{equation}
P_{-1/2+ i k}(x+1)= m(k) x^{-1/2+i k}+\overline{ m(k)}x^{-1/2-i k} + O(x^{-3/2}), \q x\to \infty,
\label{eq:La1}\end{equation}
where 
\begin{equation}
m(k )=\frac{\Gamma(i k)}  { \sqrt{2\pi}\Gamma(1/2 + i k)  }2^{i k}.
\label{eq:La2}\end{equation}

The operator $L$
is symmetric in the space $L^2 ({\Bbb R}_{+})$ on the domain $C_{0}^\infty ({\Bbb R}_{+})$, but it is not essentially self-adjoint. Since 
both functions
$P_{-1/2+ i k }(x+1)$ and $Q_{-1/2+ i k}(x+1)$   
 belong to $L^2$ in a neighborhood of the point $x=0$,  the defect indices of the operator $L$ are $(1,1)$. One of  self-adjoint extensions of $L$ from $C_{0}^\infty ({\Bbb R}_{+})$ (it will be also denoted by $L$) is defined 
on the domain ${\cal D}( L)$ consisting of functions $f(x)$ from the Sobolev class ${\HH}^2_{loc} ({\Bbb R}_{+})$ satisfying     the  boundary conditions
\begin{equation}
\exists \lim_{x\to 0}f(x) , \q f'(x)=o(x^{-1/2}), \; x\to 0,
\label{eq:BCP}\end{equation}
  (we call these boundary conditions regular); it is also required that $f\in L^2 ({\Bbb R}_{+})$ and $Lf\in L^2 ({\Bbb R}_{+})$.  Actually, the direct integration by parts shows that the operator $L$ is symmetric. Furthermore, using the appropriate Green function, we find that for all $h\in L^2 ({\Bbb R}_{+})$ the equation $(pf')'=h$ has a solution satisfying condition \e{eq:BCP}. Thus the image of the operator $L$ coincides with $ L^2 ({\Bbb R}_{+})$, and hence $L$ is self-adjoint
  (cf. \S 132, part~II, of  \cite{AhGl}).

\medskip

{\bf 3.2.}
For a study of the operator $L$, it is convenient to make a standard (see, e.g., book \cite{Titch} by E. C. Titchmarsh) change of variables. Set
\begin{equation}
t=\omega (x)= \int_{0}^x p ( y)^{-1/2}dy\q {\rm and}\q f(x) = \omega' (x)^{1/2} \tilde{f}(\omega(x))=: (F\tilde{f}) (x).
 \label{eq:C}\end{equation}
The operator $F$ is unitary in the space $L^2 ({\Bbb R}_{+})$,
and the operator $\wt{L}= F^{-1} L F$ acts by the formula
$
\wt{L}= -  d^2/dt^2  +q(\eta (t))
$ 
 where
 \[
q(x) = -16^{-1} p(x)^{-1} p' (x)^2+ 4^{-1}   p'' (x) 
\]
and $  \eta  =\omega^{-1}$ is the inverse function to  $\omega$    (so that  $x = \eta (t)$).

 In the case $p(x)= x^2+ 2x$
 we have
 \begin{equation}
 \omega (x)=  2 \ln \big(x^{1/2}+ (x+2)^{1/2}\big)- \ln 2
 \label{eq:Cx}\end{equation} 
 and hence
 \begin{equation}
\wt{L}= - \frac{d^2}{dt^2} +\tilde{q}(t)  +1/4,  
 \label{eq:C1}\end{equation}
  where $ \tilde{q}(t)= -4^{-1}\big( \eta^2 (t)+2 \eta (t)\big)^{-1}$.
  Since $ \omega (x)=  (2x)^{1/2} +O (x )$ as $x\to 0$ and $ \omega (x)= \ln (2x) +O (x^{-1})$ as $x\to \infty$, we see that
  $\eta (t)\sim t^2/2$ as $t\to 0$ and $\eta (t)\sim e^t/2$ as $t\to \infty$. It follows that $\tilde{q} (t)\sim -(4t^2)^{-1}$ as $t\to 0$ and $\tilde{q} (t)= O (e^{-t})$ as $t\to \infty$. 
 Note that the operator  $\wt{L}$ is  self-adjoint   
on the domain ${\cal D}( \wt{L})$ consisting of functions $\tilde{f}(t)$ from the Sobolev class ${\HH}^2_{loc} ({\Bbb R}_{+})$ satisfying     the boundary conditions 
\begin{equation}
\exists \lim_{t\to 0} t^{-1/2}\tilde{f}(t) ,  \q \tilde{f}'(t)-(2t)^{-1} \tilde{f}(t)=o(t^{ 1/2}), \; t\to 0,
 \label{eq:BC}\end{equation}
 and such that $ \tilde{f}\in L^2 ({\Bbb R}_{+})$, $\wt{L}\tilde{f}\in L^2 ({\Bbb R}_{+})$.

 All usual results of spectral and scattering theories can be    applied to  the operator $\wt{L}$ and then   used for   the operator $L$.  The operator $\wt{L}$  has a simple absolutely continuous spectrum coinciding with the interval $[1/4, \infty)$. It does not have eigenvalues because the equations $\wt{L}\tilde{u}=\mu \tilde{u}$, or equivalently $Lu =\mu u$, for $\mu\in{\Bbb R}$ do not have solutions from $L^2({\Bbb R}_{+})$ satisfying the regular boundary conditions at zero.
The diagonalization of the operator
 $\wt{L}$ can be constructed (see, e.g., \cite{Titch,Ya}) in the following way. Let $\tilde{u}_{k}(t ) $, $ k >0$,  be a real-valued solution of the equation 
 \begin{equation}
 \wt{L}\tilde{ u}_{k}=(k^2+1/4)\tilde{ u}_{k}
  \label{eq:C1a}\end{equation}  satisfying   boundary conditions \e{eq:BC}. It has the asymptotics
\begin{equation}
\tilde{u}_{k}(t )= m(k)e^{i k t}+ \overline{m(k)}e^{-i k t}+ o(1)
 \label{eq:C2}\end{equation} 
 as $t\to \infty$. Then the operator $\wt{U}$ defined by the equation
\begin{equation}
(\wt{U} \tilde{f})(k)= (2\pi)^{-1/2} |m(k )|^{-1} \int_{0}^\infty \tilde{u}_{k}(t ) \tilde{f}(t)dt,
 \label{eq:C3}\end{equation} 
 is unitary in the space $L^2 ({\Bbb R}_{+})$ and $(\wt{U} \wt{L}\tilde{f})(k)=(k^2+1/4)(\wt{U} \tilde{f})(k)$.

Let us now make the change of variables \e{eq:C} and set $U= F \wt{U} F^{-1}$. Note that 
\[
(2\pi)^{-1/2} |m(k )|^{-1} =\sqrt {k \tanh  \pi k  }
\]
for the function $m(k)$ defined by equation  \e{eq:La2}. It follows that
the operator $ U$ defined by the equation
\begin{equation}
( U   f )(k)= \sqrt {k \tanh \pi k  } 
\int_{0}^\infty P_{-1/2+ i k}(x+1) f(x)dx,
 \label{eq:C4}\end{equation} 
 is unitary in the space $L^2 ({\Bbb R}_{+})$ and 
 \begin{equation}
 ( U L f)(k)=(k^2+ 1/4) ( U f)(k).
  \label{eq:Int}\end{equation}

\medskip

{\bf 3.3.}   
Now we return to the Hankel operator $A$. Observe that the function $ P_{-1/2+ i k}(x+1)$ satisfies\footnote{We are obliged to choose regular boundary conditions at zero  since the function $ Q_{-1/2+ i k}(x+1)$ does not satisfy the second boundary condition \e{eq:L4}} both boundary conditions \e{eq:L4} and \e{eq:L4a}.
  It follows from Theorem~\ref{pr1} that
\begin{equation}
\int_{0}^\infty (x+y+2)^{-1} P_{-1/2+ i k}(y+1) dy = \lambda P_{-1/2+ ik}(x+1)+ \check{\lambda} Q_{-1/2+ ik}(x+1).
\label{eq:A3} \end{equation}
Considering here the limit $x\to 0$, we see that $\check{\lambda}=0$. Then we take the limit $x\to \infty$. It easily follows from  \e{eq:La1} that the left-hand side of \e{eq:A3} equals
\begin{align*}
2 \Re \Big(m(k) \int_{0}^\infty &(x+y+2)^{-1} y^{-1/2+ i k}  dy \Big)  +O(x^{-1})
\\
&=
2 \Re \Big(m(k) x^{-1/2+ i k} \int_{0}^\infty ( t+ 1)^{-1} t^{-1/2+ i k}  dt\Big) +O(x^{-1}),
\end{align*}
where we have set $y=xt$. Comparing this asymptotics with asymptotics \e{eq:La1}
of the right-hand side of \e{eq:A3}, we see that
\begin{equation}
\lambda=\int_{0}^\infty (t+1)^{-1} t^{-1/2+i k} d t= \pi (\cosh \pi k )^{-1}
\label{eq:A5} \end{equation}
and hence
\begin{equation}
\int_{0}^\infty (x+y+2)^{-1} P_{-1/2+ ik }(y+1) dy = \pi (\cosh \pi k )^{-1} P_{-1/2+ i k}(x+1). 
\label{eq:A6} \end{equation}
It yields equation \e{eq:CC}
with the operator $U$   defined by formula \e{eq:C4}. Since the operator $U$ is unitary, we have recovered the result of   F.~Mehler \cite{Me}.

\begin{proposition}\label{Me}
The Hankel operator with kernel $a(x)=(x+2)^{-1}$ has the simple absolutely continuous spectrum coinciding with the interval $[0,\pi]$. Its normalized eigenfunction corresponding to the spectral parameter $\lambda=\pi (\cosh \pi k )^{-1}$ is given by formula \e{eq:M1}.
\end{proposition}

We emphasize that   equation \e{eq:A6}
has been  obtained as a direct consequence of the commutator method, without any use of the theory of special functions.
  
 %%%%%%%%%%%%%%%%%%%%%%%%%%%
\section{ Singularity at zero}  
%%%%%%%%%%%%%%%%%%%%%%%%%%%%
 
 In the first three subsections  we study the Hankel operator   with kernel \e{eq:W1}
 and in   subs.~4  --  with kernel \e{eq:H2}. In both cases $a(x)\sim x^{-1}$ as $x\to 0$ and $a(x)$ decays exponentially as $x\to \infty$. The corresponding operator $L$ is defined by formula \e{eq:L1} where $\gamma=0$.
 
 \medskip
  
 % Recall that  the Whittaker function $ W_{-\beta, 1/2}(x )$ can be defined as a solution of equation  \e{eq:E3} for $\alpha=1/4$ such that $W_{-\beta, 1/2}( x)= e^{-x/2} x^{-\beta}(1+ O(x^{-1}))$ as $ x\to \infty$.   
 
 {\bf 4.1.}
 Note that in the case
  $\gamma=0$,   after  a change of variables 
\[
  \psi (x) = x ^{-1} \varphi(x ) 
\] 
in \e{eq:Y1},
we get again  (cf. equation \e{eq:E3}) the Schr\"odinger equation  
\begin{equation}
 -  \varphi''(x)  + (\alpha  +\beta x^{-1}- \mu  x^{-2}) \varphi(x)=0   
\label{eq:Y3} \end{equation}
with the Coulomb potential but with a non-zero orbital term. Below we set $\alpha=1/4$. 

Recall that
  the Whittaker function $ W_{-\beta, p}( x )$ can be defined as the solution of equation  \e{eq:Y3} for $  \mu=1/4 - p^2$  such that 
  \begin{equation}
  W_{-\beta, p}( x )= x^{-\beta} e^{- x /2}  (1+ O(x^{-1}))
  \label{eq:AS} \end{equation}
   as $ x\to \infty$. Of course, $  W_{-\beta, -p}( x)=  W_{-\beta, p}(x)$. In particular, the function $b(x)= W_{-\beta, 1/2}(x)$ satisfies equation \e{eq:E3}   (where $\alpha=1/4$). 
   
   As far as asymptotics as $x\to 0$ are concerned (see \S 6.8 of \cite{BE}), we note that
\begin{equation}
W_{-\beta, i k}( x)= m(k) x^{1/2+ik}+ \overline{m(k)} x^{1/2-ik}+ O (x^{3/2 }), \q k>0, \q x\to 0,
 \label{eq:C3b}\end{equation} 
where
\begin{equation}
m(k)= \Gamma (-2ik)  \Gamma^{-1} (1/2-ik+\beta).
 \label{eq:C3bm}\end{equation} 
If $p\geq 0$ and   $-1/2+p+\beta\neq   -1, -2, \ldots$, we have as $x\to 0$ 
\begin{equation}
\begin{split}
W_{-\beta, p}(x)& \sim  \Gamma (2p) \Gamma(1/2+p+\beta)^{-1} x^{1/2-p},\q p>0 ,
\\
W_{-\beta, 0}(x)& \sim -   \Gamma(1/2+ \beta) x^{1/2 } \ln x. 
\end{split}
\label{eq:om}\end{equation}
 
If   $-1/2+p+\beta=-n$ where $n=1,2,\ldots$, then taking into account formulas (6.9.4) and (6.9.36) of \cite{BE}, we can express   the Whittaker functions     in terms of the Laguerre polynomials:
  \begin{equation}
W_{-\beta, p}(x) = (-1)^{n-1} (n-1)!  e^{-x/2}x^{p+1/2} L^{2p}_{n-1} (x) .
 \label{eq:eig}\end{equation}

  If $\gamma=0$ and  $\alpha=1/4$, then
 \begin{equation}
L = - \frac{d}{dx}x^2  \frac{d}{dx} +   x ^2 /4 + \beta  x.
 \label{eq:L1b}\end{equation}
We emphasize that the coefficient $\beta $ may be arbitrary. It turns out that
 the strong degeneracy of the function $x^2$ at $x=0$ gives rise to a branch of the absolutely continuous spectrum of the operator $L$.

First, let us    
 define $L$ as a  self-adjoint  operator in the space $L^2 ({\Bbb R}_{+})$. We will check that the operator $L$   is essentially self-adjoint on the domain $C_{0}^\infty ({\Bbb R}_{+})$. Let\footnote{Note that the integral in \e{eq:C} diverges for $p(x)=x^2$ and hence the definition of the operator $F$ should be changed.}  $(F\tilde{f})(x)=x^{-1/2} \tilde{f}(\ln x)$. Then the transformation $F:  L^2({\Bbb R} )\to L^2 ({\Bbb R}_{+})$ is
  unitary and the operator $\wt{L}= F^{-1} L F$ acts by   formula \e{eq:C1}  where $ \tilde{q}(t) =    e^{2t}/4+\beta   e^t$. 
 This is already a standard Sturm-Liouville operator  in the space $\wt{\cal H}=L^2({\Bbb R} )$. 
 The potential $ \tilde{q}(t)$ tends to $0$ as $t\to -\infty$ and to $+\infty$ as $t\to+\infty$. In particular, $\wt{L}$ is   essentially self-adjoint  on   $C_{0}^\infty ({\Bbb R} )$ which implies that
$L$ is   essentially self-adjoint  on  $C_{0}^\infty ({\Bbb R}_{+} )$  in the space $ L^2({\Bbb R}_{+} ) $. Thus, a boundary condition at the point $x=0$ is unnecessary. Since 
$\tilde{q}(t)\to \infty$ as $t\to\infty$, a quantum particle can evade to $-\infty$ only. This ensures that  the spectrum of the operator $\wt{L}$ is  simple.

\medskip

{\bf 4.2.}
The expansion   over eigenfunctions of the operator $\wt{L}$ can be performed by the following standard procedure (see, e.g., \cite{Ya}, \S 5.4). As we will see later,  in addition to the simple absolutely continuous spectrum which coincides with $[1/4, \infty)$, for $\beta<-1/2$ the operator  $\wt{L}$   has a finite number of simple eigenvalues $\mu_{1},\ldots, \mu_{N}$, $N=N(\beta)$,  lying below the point $1/4$. We denote by  
$\wt{\cal H}^{(p)}$ the subspace spanned by the corresponding eigenfunctions.
Let $\tilde{u}_{k}(t ) $, $ k >0$,  be a real-valued solution of   equation 
  \e{eq:C1a}  belonging to $L^2 ({\Bbb R}_{+} )$. It has   asymptotics
 \e{eq:C2}   as $t\to -\infty$ with a function $m(k)$ which will be calculated later. Then the operator 
 $\wt{U} : \wt{\cal H} \to L^2 ({\Bbb R}_{+})$ defined by the equation (cf. \e{eq:C3})
\begin{equation}
(\wt{U} \tilde{f})(k)= (2\pi)^{-1/2} |m(k )|^{-1} \int_{-\infty }^\infty \tilde{u}_{k}(t ) \tilde{f}(t)dt,
 \label{eq:C3a}\end{equation} 
 is bounded, $\wt{U}\big|_{\wt{\cal H}^{(p)}}=0$, the mapping $\wt{U} : \wt{\cal H} \ominus \wt{\cal H}^{(p)}\to L^2 ({\Bbb R}_{+})$ is unitary and equation 
 \[
  ( \widetilde{U} \widetilde{L} \tilde{f})(k)=(k^2+ 1/4) ( \widetilde{U} \tilde{f} )(k)
  \]
  holds.
 
 The functions $u_{k}(x)= x^{-1/2} \tilde{u}_{k} (\ln x)$, $k>0$, satisfy   the equation 
\begin{equation}
 -(x^2  u_k' (x))' + 4^{-1}x^2  u_k (x)+ \beta x u_k (x)=(k^2+1/4) u_{k}(x)
 \label{eq:pps}\end{equation} 
  and can be expressed in terms of   Whittaker functions:
\begin{equation}
u_k (x)= x^{-1}W_{-\beta, i k}( x). %, \q \beta\neq -1, -2, \ldots.
 \label{eq:C3c}\end{equation}
 It follows from \e{eq:C3b} that the function $  \tilde{u}_{k} (t)=e^{t/2} u_{k}(e^t)$ has as
 $t\to-\infty$ asymptotics \e{eq:C2} with the function $m(k)$ defined by \e{eq:C3bm}.  Calculating $| m(k) |$ and making in \e{eq:C3bm} the change of variables $t=\ln x$, we find that
 the operator $ U=F \widetilde{U} F^{-1}$ is given by the equation
\[
( U   f )(k)=  \pi^{-1} \sqrt{k    \sinh 2\pi k}  | \Gamma (1/2 - ik +\beta)| \int_{0}^\infty x^{-1}W_{-\beta, i k}( x) f(x)dx.
 \] 
It  is bounded, $ {U}\big|_{ {\cal H}^{(p)}}=0$, the mapping $ {U} :  {\cal H} \ominus  {\cal H}^{(p)}\to L^2 ({\Bbb R}_{+})$ is unitary and equation \e{eq:Int} holds.
  Here $ {\cal H}^{(p)}$ is  the subspace spanned by   the  eigenfunctions  $ {\psi}_{1}, \ldots,  {\psi}_{N}$ of the operator $L$.

 Let us calculate  these functions.    The function 
 $u_{p}(x)= x^{-1} W_{-\beta,p}(x) $ for $p\geq 0$   satisfies  equation \e{eq:pps} where the role of $k^2$ is played by $-p^2$.  In view of \e{eq:AS} it belongs to $L^2$ at infinity. However, it follows from asymptotics \e{eq:om} that it does not belong to $L^2$ in a neighborhood of the point $x=0$ unless    $-1/2+p+\beta=-n$ where $n=   1,  2, \ldots$.
 Moreover, in view of \e{eq:eig} for all  $\beta=-1/2$  the function $u_{0}\not\in L^2$. Thus, if $\beta\geq -1/2$, the operator $L$ is purely absolutely continuous. 
 If $\beta<-1/2$, it also has the eigenvalues $\mu_{n}=1/4- ( |\beta| + 1/2 -n )^2$ where $n=   1, 2, \ldots$ and $n<|\beta| + 1/2$ (then $p>0$). According to
   formula \e{eq:eig} the  corresponding eigenfunctions equal
 \begin{equation}
\psi_{n} (x)=   e^{-x/2}x^{p-1/2} L^{2p}_{n-1} (x), \q p=|\beta| +1/2-n.
 \label{eq:C3n}\end{equation}
  
 \medskip

{\bf 4.3.}
Now we return to the Hankel operator $A$ with kernel  \e{eq:W1}.
It follows from \e{eq:AS} that $a(x)$ exponentially decays as $x\to\infty$,  and it follows from the first formula \e{eq:om} for $p=1/2$ that $a(x)\sim x^{-1}$ as $x\to 0$. Observe that  in view of asymptotics \e{eq:AS} and \e{eq:C3b},
function \e{eq:C3c}  satisfies  both boundary conditions \e{eq:L4} and \e{eq:L4a}.
  Hence it follows from Theorem~\ref{pr1} that
  \begin{equation}
 \int_{0}^\infty a(x+y)  y^{-1}W_{-\beta, i k}( y) dy
=  \lambda ( k) x^{-1}W_{-\beta, i k}( x)+  \check{\lambda} (k) x^{-1}M_{-\beta, i k}( x) 
 \label{eq:Sh1}\end{equation}  
 where the Whittaker function $M_{-\beta, i k}$ is the solution of equation \e{eq:E3}
  exponentially growing as $x\to\infty$. Therefore considering   the limit $x\to \infty$ in \e{eq:Sh1}, we see that necessarily $\check{\lambda}(k)=0$. Then we take the limit $x\to 0$ and use asymptotics \e{eq:C3b}. Since $a(x)\sim x^{-1}$ as $x\to 0$, we have
\begin{align*}
\int_{0}^\infty a(x+y)    y^{-1}W_{-\beta, i k}( y) dy 
= & 2 \Re \Big(m (k)  \int_{0}^\infty (x+y)^{-1}    y^{-1/2+ik} dy\Big) + O(x^{1/2})
 \\
= & 2 \lambda (k) \Re \big( m (k) x^{-1/2+ik}\big) + O(x^{1/2})
 \end{align*} 
where $ \lambda(k)$ is  again given  by formula
\e{eq:A5}. This yields equation \e{eq:Sh}.

It remains to calculate eigenvalues $\lambda_{1},\dots, \lambda_{N}$ of the operator $A$. The corresponding eigenfunctions are given by formula  \e{eq:C3n}. We proceed again from equation \e{eq:Sh1} where the role of $ik$ is played by $p= 
|\beta|+1/2 -n$. As before considering  the limit $x\to \infty$, we   see that $\check{\lambda}_{n}=0$ and hence
 \begin{equation}
 \int_{0}^\infty a (x+y)  e^{-y /2}y^{p -1/2} L^{2p}_{n-1} (y)dy=
\lambda_{n} e^{-x/2}x^{p-1/2} L^{2p}_{n-1} (x). 
 \label{eq:Sh2}\end{equation} 
   It follows from  \e{eq:W1} and \e{eq:AS} that the left-hand side here equals
  \[
\Gamma(1+\beta) x^{-\beta-1}e^{-x /2}   \int_{0}^\infty   e^{-y }y^{p-1/2} L^{2p}_{n-1} (y)dy \:\big(1+ O(x^{-1})\big),\q x\to\infty.
 \]
 Putting together formulas (2.8.46) and (10.12.33) of \cite{BE}, we see that
  \begin{equation}
 (n-1)! \int_{0}^\infty   e^{-y }y^{p-1/2} L^{2p}_{n-1} (y)dy = \Gamma (p +n-1/2).
  \label{eq:GG}\end{equation} 
  Recall also that $ L^{2p}_{n-1} (x)$ is a polynomial of degree $n-1$ with the coefficient
  $(-1)^{n-1}/(n-1)! $ at $ x^{n-1}$. Hence
  it follows from relation \e{eq:Sh2} that
\begin{equation}
 \lambda_{n}= (-1)^{n } \pi / \sin  \pi \beta, \q n=1,2,\ldots, \q n< |\beta|+1/2.
 \label{eq:Sh3}\end{equation} 

  Since eigenfunctions of the operator $L$ are orthogonal and complete, we have recovered the result of   M.   Rosenblum \cite{Ro}.

\begin{proposition}\label{Ro}
The Hankel operator $A$ with kernel \e{eq:W1} has the simple absolutely continuous spectrum coinciding with the interval $[0,\pi]$. Its normalized eigenfunction corresponding to a   point $\lambda=\pi (\cosh \pi k )^{-1}$ from the continuous spectrum   is given by the formula 
 \[
\psi_k (x)=  \pi^{-1} \sqrt{k \sinh 2\pi k } | \Gamma (1/2 - ik +\beta)|  \; x^{-1}W_{-\beta, i k}( x)   , \q k>0.
\]
Moreover, if $\beta<-1/2$, then the operator $A$ has eigenvalues \e{eq:Sh3} with the corresponding eigenfunctions defined by \e{eq:C3n}.
\end{proposition}

 \medskip

  {\bf 4.4.}
Next, we turn to the Hankel operator with singular kernel \e{eq:H2} which probably was not considered in the literature. Recall that the MacDonald function is defined by the relation
$ K_p(z)= 2^{-1} i e^{\pi i p/2} H_p^{(1)} (iz)$
 where $H_p^{(1)}$ is the Hankel function. Now the function $b(x)=x a(x)$ satisfies the Schr\"odinger equation \e{eq:E1} for the zero energy $\alpha=0$ and the coupling constant $\beta=2$. Of course, we could have taken arbitrary    $\beta>0$, but we have to exclude negative $\beta$ since in this case the function $b(x) $ grows as $x\to\infty$. 
 
 It follows from the well-known properties  of $H_p^{(1)}$  that function \e{eq:H2} has asymptotics 
    \begin{equation}
  a(x)= 4  \pi^{1/2}  x^{-3/4} e^{-\sqrt{8x}}(1+ O(x^{-1/2}))
   \label{eq:AS1} \end{equation}
   as $x\to\infty$ and $a(x)\sim x^{-1}$ as $x\to 0$.

 The corresponding operator
 \[
 L= -\frac{d}{dx}x^2\frac{d}{dx} + 2x
 \]
 can be studied quite similarly to operator \e{eq:L1b}.  For example, the operator $\widetilde{L}= F^{-1} L F$ acts by formula \e{eq:C1} where $\tilde{q}(t) =2 e^t$.
 A solution of equation \e{eq:Y1} where  $\mu= k^2 + 1/4$ belonging to $L^2$ at infinity can be expressed again  in terms of the MacDonald function 
    \[
u_k(x)=   x^{-1/2} K_{2ik} (\sqrt{8x}) . 
 \]
  According to formulas (7.2.12) and (7.2.13) of \cite{BE} we have
 \[
 u_k (x)=  m(k)  x^{-1/2+ i k } + \overline{m(k)}  x^{-1/2 - i k }+ O(x^{1/2}), \q x \to 0,
 \]
 where
 \[
 m(k)=i \pi 2^{-1+ ik} \big( \Gamma  (1+ 2 ik)  \sinh 2\pi k\big)^{-1}.
 \]
 Calculating $|m(k)|$ and using
     \e{eq:C3a}, we see that   formula \e{eq:U} now looks as
\begin{equation}
( U   f )(k)= 2  \pi^{- 1} \sqrt{k    \sinh 2\pi k}    \int_{0}^\infty 
x^{-1/2} K_{2 ik} (\sqrt{8x}) f(x)dx.
 \label{eq:C4a1}\end{equation}
 The operator $L$ does not have eigenvalues because the functions  $  x^{-1/2} K_{2p} (\sqrt{8x})$ for $p\geq 0$ do not belong to $L^2$ in a neighborhood of the point $x= 0$.
 Thus,  similarly  to subs.~4.2, we see  that the operator $ U$ defined by  formula \e{eq:C4a1} is  unitary in the space $L^2 ({\Bbb R}_{+})$ and the operator $ L$  has the simple absolutely continuous spectrum $[1/4, \infty)$.

 Theorem~\ref{pr1}  implies that 
  \begin{equation}
\int_{0}^\infty a(x+y)  y^{-1/2} K_{2 i k} (\sqrt{8y}) dy=\lambda (k) x^{-1/2} K_{2 i k} (\sqrt{8x}) + \check{\lambda} (k ) x^{-1/2} H_{2 i k}^{(2)} (i \sqrt{8x})
 \label{eq:H4}\end{equation}
 (the Hankel function  $H_{2 i k}^{(2)}(iz)$    exponentially increases as $z\to \infty$) for some constants $\lambda  (k)$ and $\check{\lambda}  (k)$. Since the integral in \e{eq:H4} (exponentially) decays as $x\to \infty$, necessarily $\check{\lambda}  (k)=0$. Comparing the asymptotics of the left- and right-hand sides of \e{eq:H4} as $x\to 0$ and using that   $a(x)\sim x^{-1}$ as $x\to 0$, we find that the constant $\lambda (k)$ is again given by formula \e{eq:A5}. Thus, similarly to the previous subsection, we obtain
 
 \begin{proposition}\label{Be}
The Hankel operator $A$ with kernel \e{eq:H2} has the simple absolutely continuous spectrum coinciding with the interval $[0,\pi]$. Its normalized eigenfunction corresponding to a spectral point $\lambda=\pi (\cosh \pi k )^{-1}$ is given by the formula 
 \[
\psi_k (x)= 2  \pi^{- 1} \sqrt{k    \sinh 2\pi k}    x^{-1/2} K_{2 ik} (\sqrt{8x})   , \q k>0.
\]
\end{proposition}

As a by-product of our considerations, we obtain the equation
 \[
x^{1/2}  \int_{0}^\infty (x+y)^{-1/2} K_1 (\sqrt{ x+y  }) y^{-1/2} K_{2i k} (\sqrt{ y}) dy = \pi (\cosh   \pi k  )^{-1}   K_{2 i k} (\sqrt{x}).
\]
We have not found this equation in the literature on special functions. Note, however, that it can   formally  be deduced from the Shanker equation \e{eq:Sh} if one uses the relation (formula (6.9.19) of \cite{BE})
\[
\lim_{\beta \to \infty}   \Gamma (\beta +1)  W_{-\beta, m} (x/\beta)= 2 x^{1/2} K_{2m}(2 x^{1/2}).
\]

\medskip
 
 {\bf 4.5.}
 The Carleman operator $A$ trivially fits into the scheme exposed above. Now the operator $A$ commutes with   operator \e{eq:L1} for $\alpha=\beta=\gamma=0$. This operator has the absolutely continuous spectrum of multiplicity $2$ coinciding with $[1/4,\infty)$. It has eigenfunctions $x^{-1/2+ ik}$ for all $k\in {\Bbb R}$ which are also eigenfunctions of the operator $A$. The relation between the spectral parameters $\lambda$ and $k$ is again given by formula  \e{eq:A5} so that
 the operator $A$  has the absolutely continuous spectrum of multiplicity $2$ coinciding with $[0,\pi]$.

 %%%%%%%%%%%%%%%%%%%%%%%%%%%
\section{Regular kernels}  
%%%%%%%%%%%%%%%%%%%%%%%%%%%%

 {\bf 5.1.}
  Let us here consider kernels $a(x)$ which decay rapidly as $x\to \infty$ and have finite limits as $x\to 0$. We set   $\gamma=2$ and distinguish the cases $\alpha> 0$, $\beta$ is arbitrary and $\alpha=0$, $\beta>0$.  Let $\alpha=1/4$ and $\beta=2 $ in the first and second cases, respectively. If $\alpha=1/4$, then
  the solution of equation \e{eq:E3}    is given    (see subs.~4.1)   by the formula $b(r)= W_{-\beta, 1/2}(r)$ where $W_{-\beta, 1/2}$ is the Whittaker function.  If $\alpha=0$ and $\beta=2 $, then
  the solution of   \e{eq:E3}   equals $b(r)= r^{1/2} K_{1}(\sqrt{8r})$  where  $K_{1} $ is the MacDonald function  (see subs.~4.4).  The corresponding functions \e{eq:E2}  decay exponentially at infinity  and   have finite limits   as $x\to 0$.  It follows that the operators $A$ are compact. 

Let  the function $\omega(x)$ be defined by formula \e{eq:Cx}, $\eta=\omega^{-1}$    and
\begin{equation}
\tilde{q}(t)= -4^{-1}\big( \eta^2 (t)+2 \eta (t)\big)^{-1}+\alpha  \eta^2(t)+\beta \eta (t).
\label{eq:qq} \end{equation}
Since  $\eta(t)\sim e^t/2$, the potential $\tilde{q}(t)\to +\infty$  as $ t\to\infty$. It follows that
   the   operators $\widetilde{L}$  and hence $L$ have now discrete spectra. We point out that these operators are again defined by formulas \e{eq:C1} and \e{eq:L1}  on functions satisfying boundary conditions \e{eq:BC} and \e{eq:BCP}, respectively.

  Theorem~\ref{pr1} implies that the  operators $A$ and $L$  have common eigenfunctions. Apparently,  eigenfunctions of the  operator  $L$ cannot be expressed in terms of standard special functions. However, in  their  terms we can calculate eigenvalues $\lambda_{1}, \lambda_{2},\ldots$ of the operator $A$. Indeed, suppose  that $\psi_{\mu}\in{\cal D}(L)$ and $L\psi_{\mu}=\mu \psi_{\mu}$.
   Then the function $\tilde{\psi} _{\mu}(t)= ({\sf F}^{-1} \psi_{\mu})(t)$ satisfies   the equation 
  \[
  \tilde{\psi} _{\mu}''(t)+ \tilde{q}(t) \tilde{\psi} _{\mu}(t) =(\mu-1/4)\tilde{\psi} _{\mu}(t)
  \]
where  $\tilde{q}(t)$ is function \e{eq:qq}. For a suitable normalization, asymptotics of $ \tilde{\psi} _{\mu}(t)$ as $t\to \infty$  is given (see, e.g., book  \cite{Olver}) by the semiclassical formula
  \[
  \tilde{\psi} _{\mu}(t)\sim \tilde{q}(t)^{-1/4}\exp \Big(-\int_{0}^t \tilde{q}(s)^{1/2} ds\Big) .
  \]
  It follows that $ \tilde{\psi} _{\mu}(t)\sim e^{-t/2} \exp(-4^{-1}e^t-\beta t)$ in the first case and $ \tilde{\psi} _{\mu}(t)\sim e^{-t/4} \exp(-2^{3/2}e^{t/2})$  in the second case. Returning to the eigenfuctions $\psi_{\mu}(x)$, we find that
 \begin{equation}
  \psi_{\mu}(x)\sim x^{-1-\beta}e^{-x/2}\q {\rm and}\q 
   \psi_{\mu}(x)\sim x^{-3/4}e^{-\sqrt{8x}  },\q x\to\infty, 
 \label{eq:ps} \end{equation}
  in the first and second cases, respectively. 
   
   On the other hand, using asymptotics \e{eq:AS} and \e{eq:AS1} for function \e{eq:H2}, we see that
   \[
   (A\psi_{\mu} )(x)\sim x^{-1-\beta}e^{-x/2}\int_{0}^\infty e^{-y/2} \psi_{\mu} (y)dy
   \]
   and
     \[
   (A\psi_{\mu} )(x)\sim \sqrt{2\pi} x^{-3/4}e^{-\sqrt{8x}}\int_{0}^\infty   \psi_{\mu} (y)dy
   \]
   as $x\to\infty$  in the first and second cases, respectively. Comparing these relations with relations  \e{eq:ps} and using the equation $A\psi_{\mu}=\lambda_{\mu} \psi_{\mu}$, we get expressions for eigenvalues of the operators $A$:
   \begin{equation}
\lambda_{\mu} = \int_{0}^\infty e^{-y/2} \psi_{\mu} (y)dy \q {\rm and}\q 
\lambda_{\mu} = \sqrt{2\pi} \int_{0}^\infty   \psi_{\mu} (y)dy
 \label{eq:pss} \end{equation}
 in the first and second cases, respectively. 
   
  Thus we have obtained the following results.

\begin{proposition}\label{MD}
The Hankel operator $A$ with kernel 
$$a(x)=(x+2)^{-1}  W_{-\beta, 1/2}(x+2)$$
 and the differential operator  \e{eq:L1} for $\gamma=2$ and $\alpha=1/4$  have common eigenfunctions. 
  If $\psi_{\mu}\in{\cal D}(L)$, $L\psi_{\mu}=\mu \psi_{\mu}$ and $\psi_{\mu} (x)$  has the first   asymptotics \e{eq:ps} as $x\to\infty$, then $A\psi_{\mu}=\lambda_{\mu} \psi_{\mu}$ where $\lambda_{\mu}$ is determined  by the first   formula \e{eq:pss}.
\end{proposition}

\begin{proposition}\label{MDb}
The   Hankel operator $A$ with kernel 
 $$a(x)=(x+2)^{-1/2}    K_{1}(\sqrt{8(x+2)})$$
  and the differential operator  \e{eq:L1} for $\gamma=2$, $\alpha=0$  and $\beta=2$ have common eigenfunctions. 
  If $\psi_{\mu}\in{\cal D}(L)$, $L\psi_{\mu}=\mu \psi_{\mu}$ and $\psi_{\mu} (x)$  has the  second  asymptotics \e{eq:ps} as $x\to\infty$, then $A\psi_{\mu}=\lambda_{\mu} \psi_{\mu}$ where $\lambda_{\mu}$ is determined  by the first  second  formula \e{eq:pss}.
\end{proposition}

\medskip
 
 {\bf 5.2.}
 Finally, we consider kernel \e{eq:W1} for exceptional values $\beta=-l$ where  $l=1, 2, \ldots$. To be more precise, we now set  
   \begin{equation}
a(x)= (-1)^{l-1} (l-1)!^{-1}  x^{-1} W_{l, 1/2}(x )= e^{-x/2} L^1_{l-1}(x), \q l =1, 2, \ldots,
\label{eq:tr}\end{equation}
(here we   have    taken formula \e{eq:eig} into account).
The Hankel operator $A$ with this kernel has   rank $l$.  Here we show how this simple example fits into the scheme exposed above.

The spectral analysis of the corresponding operator \e{eq:L1b} remains the same as in subs.~4.2. In addition to the absolutely continuous spectrum $[1/4,\infty)$, the operator $L$ has eigenvalues $\mu_{n}= 1/4-(l+ 1/2-n)^2 $ where $n=1,\ldots,l$.  

However, instead of the absolutely continuous spectrum, the operator $A$ has the zero eigenvalue of infinite multiplicity. Indeed, as in subs.~4.3,
Theorem~\ref{pr1} yields    equation \e{eq:Sh1} where again $\check{\lambda}(k)=0$. Observe that $a(x)$ 
 and hence in view of \e{eq:C3b} the integral in the left-hand side of \e{eq:Sh1} have finite limits as $x\to 0$. Therefore it follows from \e{eq:C3b}   that necessarily $\lambda(k)=0$  for all $k>0$. Hence the kernel of the operator $A$ is spanned by the functions $x^{-1}W_{l, i k}( x)$, $k>0$. 
 
Eigenfunctions $\psi_{n}(x)$ corresponding to non-zero eigenvalues $\lambda_{n}$ of the operator $A$ are defined by formula \e{eq:C3n} where $p=l+ 1/2-n$, $n<l+1/2$,  and  $\lambda_{n}$ can be found from equation \e{eq:Sh2}: 
\begin{equation}
   \int_{0}^\infty     L^{1}_{l -1} (x+y)  e^{-y  }y^{p-1/2} L^{2p}_{n-1} (y)dy=
\lambda_{n}  x^{p -1/2} L^{2p}_{n-1} (x). 
 \label{eq:LG}\end{equation} 
 Recall that $L_{p}^\alpha(x)$ is a polynomial of degree $p$ with the coefficient $(-1)^p /p!$ at $x^p$. Comparing   coefficients at the highest power $x^{l-1}$ in the left- and right-hand sides of \e{eq:LG} and taking into account formula \e{eq:GG}, we find that
\[
 \lambda_{n}=(-1)^{n-l} \frac{(n-1)!}{(l-1)!} \int_{0}^\infty    e^{-y  }y^{p-1/2}   L^{2p}_{n-1} (y)dy=
 (-1)^{n-l}  .
  \]
  
  Thus, we have obtained the following result.

\begin{proposition}\label{Ro1}
The Hankel operator $A$ with kernel \e{eq:tr} has    rank $l$. Its non-zero eigenvalues
are given by the formula $\lambda_{n}= 
 (-1)^{n-l}$ where $n=1,\ldots,l$, and the corresponding eigenfunctions $\psi_{n} (x)$ are defined by equality \e{eq:C3n} where $p=l +1/2 -n$.
\end{proposition}

\end{document}